\documentclass{elsart}
\pdfoutput=1

\usepackage{graphicx}
\usepackage{ifpdf}
\ifpdf
\fi

\usepackage{amssymb, amsmath}
\usepackage{color}

\newcommand{\symm}{\bigtriangleup}

\begin{document}

\begin{frontmatter}



\title{Stability results for uniquely determined sets from two directions in discrete tomography}


\author{Birgit van Dalen}

\address{Mathematisch Instituut, Universiteit Leiden, Niels Bohrweg 1, 2333 CA Leiden, The Netherlands}

\ead{dalen@math.leidenuniv.nl}

\begin{abstract}
In this paper we prove several new stability results for the reconstruction of binary images from two projections. We consider an original image that is uniquely determined by its projections and possible reconstructions from slightly different projections. We show that for a given difference in the projections, the reconstruction can only be disjoint from the original image if the size of the image is not too large. We also prove an upper bound for the size of the image given the error in the projections and the size of the intersection between the image and the reconstruction.
\end{abstract}

\begin{keyword}

Discrete tomography \sep Stability \sep Image reconstruction \sep Binary image


\end{keyword}
\end{frontmatter}

\section{Introduction}

Discrete tomography is concerned with problems such as reconstructing binary images on a lattice from given projections in lattice directions \cite{boek}. Each point of a binary image has a value equal to zero or one. The line sum of a line through the image is the sum of the values of the points on this line. The projection of the image in a certain lattice direction consists of all the line sums of the lines through the image in this direction.

Several problems related to the reconstruction of binary images from two or more projections have been described in the literature \cite{boek, boeknieuw}. Already in 1957, Ryser gave an algorithm to reconstruct binary images from their horizontal and vertical projections and characterised the set of projections that correspond to a unique binary image \cite{ryser}. For any set of directions, it is possible to construct images that are not uniquely determined by their projections in those directions \cite[Theorem 4.3.1]{boek}. The problem of deciding whether an image is uniquely determined by its projections and the problem of reconstructing it are NP-hard for any set of more than two directions \cite{gardner}.

Aside from various interesting theoretical problems, discrete tomography also has applications in a wide range of fields. The most important are electron microscopy \cite{gold} and medical imaging \cite{medicaldt, medical}, but there are also applications in nuclear science \cite{radio2, radio} and various other fields \cite{propellor, building}.

An interesting problem in discrete tomography is the stability of reconstructions. Even if an image is uniquely determined by its projections, a very small error in the projections may lead to a completely different reconstruction \cite{alpersthesis, alpersinstabiel}. Alpers et al. \cite{alpersthesis, alpersartikel} showed that in the case of two directions a total error of at most 2 in the projections can only cause a small difference in the reconstruction. They also proved a lower bound on the error if the reconstruction is disjoint from the original image.

In this paper we improve this bound, and we resolve the open problem of stability with a projection error greater than 2.

\section{Notation and statement of the problems}\label{notation}

Let $F_1$ and $F_2$ be two finite subsets of $\mathbb{Z}^2$ with characteristic functions $\chi_1$ and $\chi_2$. (That is, $\chi_h(x,y) = 1$ if and only if $(x,y) \in F_h$, $h \in \{1,2\}$.) For $i \in \mathbb{Z}$, we define \emph{row} $i$ as the set $\{(x,y) \in \mathbb{Z}^2: x = i\}$. We call $i$ the index of the row. For $j \in \mathbb{Z}$, we define \emph{column} $j$ as the set $\{(x,y) \in \mathbb{Z}^2: y = j\}$. We call $j$ the index of the column. Following matrix notation, we use row numbers that increase when going downwards and column numbers that increase when going to the right.

The \emph{row sum} $r_i^{(h)}$ is the number of elements of $F_h$ in row $i$, that is $r_i^{(h)} = \sum_{j \in \mathbb{Z}} \chi_h(i,j)$. The \emph{column sum} $c_j^{(h)}$ of $F_h$ is the number of elements of $F_h$ in column $j$, that is $c_j^{(h)} = \sum_{i \in \mathbb{Z}} \chi_h(i,j)$. We refer to both row and column sums as the \emph{line sums} of $F_h$.

\begin{figure}
\begin{center}
\includegraphics{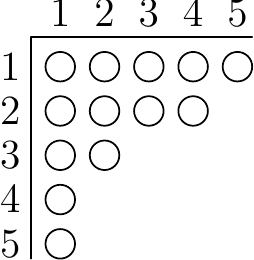}
\end{center}
\caption{A uniquely determined set with the assumed row and column ordering.}
\label{plaatjeordening}
\end{figure}

Throughout this paper, we assume that $F_1$ is uniquely determined by its row and column sums. Such sets were studied by, among others, Ryser \cite{ryser} and Wang \cite{wang}. Let $a$ be the number of rows and $b$ the number of columns that contain elements of $F_1$. We renumber the rows and columns such that we have
\[
r_1^{(1)} \geq r_2^{(1)} \geq \ldots \geq r_a^{(1)} > 0,
\]
\[
c_1^{(1)} \geq c_2^{(1)} \geq \ldots \geq c_b^{(1)} > 0,
\]
and such that all elements of $F_2$ are contained in rows and columns with positive indices. By \cite[Theorem 2.3]{wang} we have the following property of $F_1$ (see Figure \ref{plaatjeordening}):
\begin{itemize}
\item in row $i$ the elements of $F_1$ are precisely the points $(i,1)$, $(i,2)$, \ldots, $(i,r_i^{(1)})$,
\item in column $j$ the elements of $F_1$ are precisely the points $(1, j)$, $(2,j)$, \ldots, $(c_j^{(1)}, j)$.
\end{itemize}
We will refer to this property as the \emph{triangular shape} of $F_1$.

Everywhere except in Section \ref{neq} we assume that $|F_1| = |F_2|$. Note that we do not assume $F_2$ to be uniquely determined.

As $F_1$ and $F_2$ are different and $F_1$ is uniquely determined by its line sums, $F_2$ cannot have exactly the same line sums as $F_1$. Define the \emph{difference} or \emph{error in the line sums} as
\[
\sum_{j \geq 1} |c_j^{(1)} - c_j^{(2)}| + \sum_{i
\geq 1} |r_i^{(1)} - r_i^{(2)}|.
\]
As in general $|t-s| \equiv t+s$ mod 2, the above expression is congruent to
\[
\sum_{j \geq 1} \left(c_j^{(1)} + c_j^{(2)}\right) + \sum_{i \geq 1} \left( r_i^{(1)} + r_i^{(2)} \right) \equiv 2|F_1| + 2|F_2| \equiv 0 \mod 2,
\]
hence the error in the line sums is always even. We will denote it by $2 \alpha$, where $\alpha$ is a positive integer.

For notational convenience, we will often write $p$ for $|F_1 \cap F_2|$.

We consider two problems concerning stability.

\begin{prob}\label{problem1}
Suppose $F_1 \cap F_2 = \emptyset$. How large can $|F_1|$ be in terms of $\alpha$?
\end{prob}

Alpers et al. \cite[Theorem 29]{alpersartikel} proved that $|F_1| \leq \alpha^2$. They also showed that there is no constant $c$ such that $|F_1| \leq c \alpha$ for all $F_1$ and $F_2$. In Section \ref{solution1} of this paper we will prove the new bound $|F_1| \leq \alpha ( 1+ \log \alpha)$ and show that this bound is asymptotically sharp.

\begin{prob}\label{problem2}
How small can $|F_1 \cap F_2|$ be in terms of $|F_1|$ and $\alpha$, or, equivalently, how large can $|F_1|$ be in terms of $|F_1 \cap F_2|$ and $\alpha$?
\end{prob}

Alpers (\cite[Theorem 5.1.18]{alpersthesis}) showed in the case $\alpha = 1$ that
\[
|F_1 \cap F_2| \geq |F_1| + \tfrac{1}{2} - \sqrt{2|F_1|+\tfrac{1}{4}}.
\]
This bound is sharp: if $|F_1| = \frac{1}{2}n(n+1)$ for some positive integer $n$, then there exists an example for which equality holds. A similar result is stated in \cite[Theorem 19]{alpersartikel}.

While \cite{alpersthesis,alpersartikel} only deal with the case $\alpha =1$, we will give stability results for general $\alpha$. In Section \ref{solution2} we will give two different upper bounds for $|F_1|$. The bounds have different asymptotic behaviour. Writing $p$ for $|F_1 \cap F_2|$, the second bound reduces to
\[
|F_1| \leq p+1 + \sqrt{2p+1}
\]
in case $\alpha = 1$, which is equivalent to
\[
p \geq |F_1| - \sqrt{2|F_1|}.
\]
Hence the second new bound can be viewed as a generalisation of Alpers' bound. The first new bound is different and better in the case that $\alpha$ is very large.

In Section \ref{neq} we will generalise the results to the case $|F_1| \neq |F_2|$.

\section{Staircases}

Alpers introduced the notion of a staircase to characterise $F_1 \symm F_2$ in the case $\alpha = 1$. We will use a slightly different definition and then show that for general $\alpha$ the symmetric difference $F_1 \symm F_2$ consists of $\alpha$ staircases.

\begin{defn}
A set of points $(p_1, p_2, \ldots, p_n)$ in $\mathbb{Z}^2$ is called a
\emph{staircase} if the following two conditions are satisfied:
\begin{itemize}
\item for each $i$ with $1 \leq i \leq n-1$ one of the points
$p_i$ and $p_{i+1}$ is an element of $F_1 \backslash F_2$ and the
other is an element of $F_2 \backslash F_1$;
\item either for all $i$ the points $p_{2i}$ and $p_{2i+1}$ are in
the same column and the points $p_{2i+1}$ and $p_{2i+2}$ are in
the same row, or for all $i$ the points $p_{2i}$ and $p_{2i+1}$
are in the same row and the points $p_{2i+1}$ and $p_{2i+2}$ are
in the same column.
\end{itemize}
\end{defn}

This definition is different from \cite{alpersthesis,alpersartikel} in the following way. Firstly, the number of points does not need to be even. Secondly, the points $p_1$ and $p_n$ can both be either in $F_1 \backslash F_2$ or in $F_2 \backslash F_1$. So this definition is slightly more general than the one used in \cite{alpersthesis,alpersartikel} for the case $\alpha =1$.

\begin{figure}
\begin{center}
\includegraphics{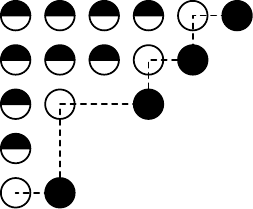}
\end{center}
\caption{A staircase. The set $F_1$ consists of the white and the black-and-white points, while $F_2$ consists of the black and the black-and-white points. The staircase is indicated by the dashed line segments.}
\label{plaatjetrap}
\end{figure}

Consider a point $p_i \in F_1 \backslash F_2$ of a staircase
$(p_1, p_2, \ldots, p_n)$. Assume $p_{i-1}$ is in the same column
as $p_i$ and $p_{i+1}$ is in the same row as $p_i$. Because of the
triangular shape of $F_1$, the row index of $p_{i-1}$ must be
larger than the row index of $p_i$, and the column index of
$p_{i+1}$ must be larger than the column index of $p_i$.
Therefore, the staircase looks like a real-world staircase (see Figure \ref{plaatjetrap}). From now on,
we assume for all staircases that $p_1$ is the point with the
largest row index and the smallest column index, while $p_n$ is
the point with the smallest row index and the largest column
index. We say that the staircase \emph{begins} with $p_1$ and \emph{ends} with $p_n$.

\begin{lem}\label{staircases}
Let $F_1$ and $F_2$ be finite subsets of $\mathbb{Z}^2$ such that
\begin{itemize}
\item $F_1$ is uniquely determined by its row and column sums, and
\item $|F_1| = |F_2|$.
\end{itemize}
Let $\alpha$ be defined as in Section \ref{notation}. Then the set $F_1 \symm F_2$ is the disjoint union of $\alpha$
staircases.
\end{lem}

\begin{pf}
We will construct the staircases one by one and delete them from
$F_1 \symm F_2$. For a subset $A$ of $F_1 \symm F_2$, define
\begin{eqnarray*}
\rho_i(A) &=& |\{ j \in \mathbb{Z}: (i,j) \in A \cap F_1 \}| - |\{
j \in \mathbb{Z}: (i,j) \in A \cap F_2 \}|, \quad i \in \mathbb{Z}, \\
\sigma_j(A) & = & |\{ i \in \mathbb{Z}: (i,j) \in A \cap F_1 \}| -
|\{ i \in \mathbb{Z}: (i,j) \in A \cap F_2 \}|, \quad j \in \mathbb{Z}, \\
\tau(A) &=& \sum_{i} |\rho_i(A)| + \sum_j |\sigma_j(A)|.
\end{eqnarray*}
We have $2\alpha = \tau(F_1 \symm F_2)$.

Assume that the rows and columns are ordered as in Section \ref{notation}. Because of the triangular shape of $F_1$, for any point
$(i,j) \in F_1 \backslash F_2$ and any point $(k,l) \in F_2
\backslash F_1$ we then have $k > i$ or $l > j$.

Suppose we have deleted some staircases and are now left with a
non-empty subset $A$ of $F_1 \symm F_2$. Let $(p_1, p_2, \ldots,
p_n)$ be a staircase of maximal length that is contained in $A$. Let $(x_1, y_1)$ and
$(x_n,y_n)$ be the coordinates of the points $p_1$ and $p_n$
respectively. Each of those two points can be either in $A \cap
F_1$ or in $A \cap F_2$, so there are four different cases. (If $n=1$, so $p_1$ and $p_n$ are the same point, then there are only two cases.) We consider two cases; the other two are similar.

First suppose $p_1 \in A \cap
F_1$ and $p_n \in A \cap F_2$. If $(x,y_1)$ is a point of $A \cap
F_2$ in the same column as $p_1$, then $x > x_1$, so we can extend the staircase by adding
this point. That contradicts the maximal length of the staircase. So there are no points of $A \cap F_2$
in column $y_1$. Therefore $\sigma_{y_1}(A) > 0$.

Similarly, since $p_n \in A \cap F_2$, there are no points of $A \cap F_1$
in the same column as $p_n$. Therefore $\sigma_{y_n}(A)
< 0$.

All rows and all columns that contain points of the staircase,
except columns $y_1$ and $y_n$, contain exactly two points of the
staircase, one in $A \cap F_1$ and one in $A \cap F_2$. Let $A' =
A \backslash \{ p_1, p_2, \ldots, p_n \}$. Then $\rho_i(A') =
\rho_i(A)$ for all $i$, and $\sigma_j(A') = \sigma_j(A)$ for all
$j \neq y_1, y_n$. Furthermore, $\sigma_{y_1}(A') =
\sigma_{y_1}(A) - 1$ and $\sigma_{y_n}(A') = \sigma_{y_n}(A) + 1$.
Since $\sigma_{y_1}(A) > 0$ and $\sigma_{y_n}(A) < 0$, this gives
$\tau(A') = \tau(A) - 2$.

Now consider the case $p_1 \in A \cap F_1$ and $p_n \in A \cap F_1$. As above, we have $\sigma_{y_1}(A) > 0$. Suppose $(x_n,y)$ is a point of $A \cap F_2$ in the same row as $p_n$. Then $y>y_n$, so we can extend the staircase by adding this point. That contradicts the maximal length of the staircase. So there are no points of $A \cap F_2$ in row $x_n$. Therefore $\rho_{x_n}(A) > 0$.

All rows and all columns that contain points of the staircase, except column $y_1$ and row $x_n$, contain exactly two points of the staircase, one in $A \cap F_1$ and one in $A \cap F_2$. Let $A' = A \backslash \{ p_1, p_2, \ldots, p_n \}$. Then $\rho_i(A') = \rho_i(A)$ for all $i \neq x_n$, and $\sigma_j(A') = \sigma_j(A)$ for all $j \neq y_1$. Furthermore, $\sigma_{y_1}(A') = \sigma_{y_1}(A) - 1$ and $\rho_{x_n}(A') = \rho_{x_n}(A) - 1$. Since $\sigma_{y_1}(A) > 0$ and $\rho_{x_n}(A) > 0$, this gives $\tau(A') = \tau(A) - 2$.

We can continue deleting staircases in this way until all points of $F_1 \symm
F_2$ have been deleted. Since $\tau(A) \geq 0$ for all subsets $A
\subset F_1 \symm F_2$, this must happen after deleting exactly
$\alpha$ staircases. \hfill $\square$
\end{pf}

\begin{rem}\label{endpoints}
Some remarks about the above lemma and its proof.
\begin{itemize}
\item[(i)] The $\alpha$ staircases from the previous lemma have $2
\alpha$ endpoints in total (where we count the same point twice in case of
a staircase consisting of one point). Each endpoint contributes a difference of 1 to the line sums in one row or column. Since all these differences must add up to $2\alpha$, they cannot cancel each other.
\item[(ii)] A staircase consisting of more than one point can be split into two or more staircases. So it may be possible to write $F_1 \symm F_2$ as the disjoint union of more than $\alpha$ staircases. However, in that case some of the contributions of the endpoints of staircases to the difference in the line sums cancel each other. On the other hand, it is impossible to decompose $F_1 \symm F_2$ into fewer than $\alpha$ staircases.
\item[(iii)] The endpoints of a staircase can be in $F_1
\backslash F_2$ or $F_2 \backslash F_1$. For a staircase $T$ of which the two
endpoints are in different sets, we have $|T \cap F_1| = |T \cap F_2|$. For a
staircase $T$ of which the two endpoints are in the same set, we have $|T \cap F_1| = 1 + |T \cap F_2|$ or $|T \cap F_2| = 1 + |T \cap F_1|$. Since $|F_1 \backslash F_2| = | F_2
\backslash F_1|$, the number of staircases with two endpoints in
$F_1 \backslash F_2$ must be equal to the number of staircases
with two endpoints in $F_2 \backslash F_1$. This implies that of the $2 \alpha$ endpoints, exactly $\alpha$ are in the set $F_1 \backslash F_2$ and $\alpha$ are in the set $F_2 \backslash F_1$.
\end{itemize}
\end{rem}

Consider a decomposition of $F_1 \symm F_2$ as in the proof of Lemma \ref{staircases}. We will now show that for our purposes we may assume that all these
staircases begin with a point $p_1 \in F_1 \backslash F_2$ and end with
a point $p_n \in F_2 \backslash F_1$.

Suppose there is a staircase beginning with a point $(x,y) \in F_2
\backslash F_1$. Then there also exists a staircase ending with a
point $(x',y') \in F_1 \backslash F_2$: otherwise more than half of the $2 \alpha$ endpoints would be in $F_2 \backslash F_1$, which is a contradiction to Remark \ref{endpoints}(iii). Because of Remark \ref{endpoints}(i) we must have $r_x^{(1)} < r_x^{(2)}$
and $r_{x'}^{(1)} > r_{x'}^{(2)}$.

Let $y''$ be such that $(x', y'') \not\in F_1 \cup F_2$. Delete the point $(x,y)$ from $F_2$ and add the point $(x', y'')$ to $F_2$. Then
$r_x^{(2)}$ decreases by 1 and $r_{x'}^{(2)}$ increases by 1, so the
difference in the row sums decreases by 2. Meanwhile, the
difference in the column sums increases by at most 2. So
$\alpha$ does not increase, while $F_1$, $|F_2|$ and $|F_1 \symm
F_2|$ do not change. So the new situation is just as good or better than the old one. The staircase that began with
$(x,y)$ in the old situation now begins with a point of $F_1
\backslash F_2$. The point that we added becomes the new endpoint
of the staircase that previously ended with $(x', y')$.

Therefore, in our investigations we may assume that all staircases begin
with a point of $F_1 \backslash F_2$ and end with a point of $F_2
\backslash F_1$. This is an important assumption that we will use in the proofs throughout the paper. An immediate consequence of it is that $r_i^{(1)}
= r_i^{(2)}$ for all $i$. The only difference between corresponding line sums occurs
in the columns.

\section{A new bound for the disjoint case}\label{solution1}

Using the concept of staircases, we can prove a new bound for Problem \ref{problem1}.

\begin{thm}\label{disjoint} Let $F_1$ and $F_2$ be finite subsets of $\mathbb{Z}^2$ such that
\begin{itemize}
\item $F_1$ is uniquely determined by its row and column sums,
\item $|F_1| = |F_2|$, and
\item $F_1 \cap F_2 = \emptyset$.
\end{itemize}
Let $\alpha$ be defined as in Section \ref{notation}. Then
\[
|F_1| \leq \sum_{i=1}^\alpha \left\lfloor \frac{\alpha}{i} \right\rfloor.
\]
\end{thm}

\begin{pf}
Assume that the rows and columns are ordered as in Section \ref{notation}. Let $a$ be the number of rows and $b$ the number of columns that contain elements of $F_1$. Let $(k,l) \in F_1$. Then all the points in the rectangle $\{ (i,j): 1 \leq i \leq k, 1 \leq j \leq l \}$ are elements of $F_1$. Since $F_1$ and $F_2$ are disjoint, none of the points in this rectangle is an element of $F_2$, and all the points belong to $F_1 \symm F_2$. So all of the $kl$ points must belong to different staircases, which implies $\alpha \geq kl$. For all $i$ with $1 \leq i \leq a$ we have $(i, r_i^{(1)}) \in F_1$, hence $r_i^{(1)} \leq \frac{\alpha}{i}$. Since $r_i^{(1)}$ must be an integer, we have
\[
|F_1| = \sum_{i=1}^a r_i^{(1)} \leq \sum_{i=1}^a \left\lfloor \frac{\alpha}{i} \right\rfloor.
\]
Since $(a,1) \in F_1$, we have $a \leq \alpha$, so
\[
|F_1| \leq \sum_{i=1}^\alpha \left\lfloor \frac{\alpha}{i} \right\rfloor.
\] \hfill $\square$
\end{pf}

\begin{cor}\label{disjunctegrenslog}
Let $F_1$, $F_2$ and $\alpha$ be defined as in Theorem \ref{disjoint}. Then
\[
|F_1| \leq \alpha (1 +\log \alpha).
\]
\end{cor}

\begin{pf}
We have
\[
|F_1| \leq \sum_{i=1}^\alpha \left\lfloor \frac{\alpha}{i} \right\rfloor \leq \alpha \sum_{i=1}^\alpha \frac{1}{i} \leq \alpha \left(1 + \int_1^\alpha \frac{1}{x}dx \right) = \alpha \left( 1 + \log \alpha \right).
\] \hfill $\square$
\end{pf}

The following example shows that the upper bound cannot even be improved by a factor $\frac{1}{2 \log 2} \approx 0.72$.

\begin{figure}
\begin{center}
\includegraphics{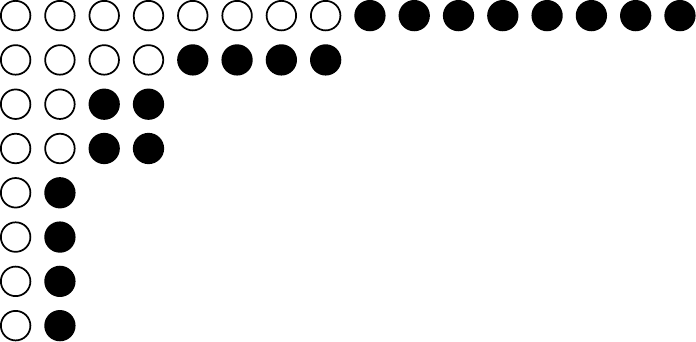}
\end{center}
\caption{The construction from Example \ref{example1} with $m=3$.}
\label{plaatjevoorbeeld1}
\end{figure}

\begin{exmp}\label{example1} (taken from \cite{alpersthesis}) \end{exmp} Let $m \geq 1$ be an integer. We construct sets $F_1$ and $F_2$ as follows (see also Figure \ref{plaatjevoorbeeld1}).
\begin{itemize}
\item Row 1: \begin{itemize}
\item $(1,j) \in F_1$ for $1 \leq j \leq 2^m$,
\item $(1,j) \in F_2$ for $2^m + 1 \leq j \leq 2^{m+1}$.
\end{itemize}
\item Let $0 \leq l \leq m-1$. Row $i$, where $2^l + 1 \leq i \leq 2^{l+1}$:
\begin{itemize}
\item $(i,j) \in F_1$ for $1 \leq j \leq 2^{m-l-1}$,
\item $(i,j) \in F_2$ for $2^{m-l-1}+1 \leq j \leq 2^{m-l}$.
\end{itemize}
\end{itemize}

The construction is almost completely symmetrical: if $(i,j) \in F_1$, then $(j,i) \in F_1$; and if $(i,j) \in F_2$ with $i > 1$, then $(j,i) \in F_2$. Since it is clear from the construction that each row contains exactly as many points of $F_1$ as points of $F_2$, we conclude that each column $j$ with $2 \leq j \leq 2^m$ contains exactly as many points of $F_1$ as points of $F_2$ as well. The only difference in the line sums occurs in the first column (which has $2^m$ points of $F_1$ and none of $F_2$) and in columns $2^m+1$ up to $2^{m+1}$ (each of which contains one point of $F_2$ and none of $F_1$). So we have
\[
\alpha = 2^m.
\]
Furthermore,
\[
|F_1| = 2^m + \sum_{l=0}^{m-1} 2^l 2^{m-l-1} = 2^m + m2^{m-1}.
\]
Hence for this family of examples it holds that
\[
|F_1| = \alpha + \frac{1}{2} \alpha \log_2 \alpha,
\]
which is very close to the bound we proved in Corollary \ref{disjunctegrenslog}.

\section{Two bounds for general $\alpha$}\label{solution2}

In case $F_1$ and $F_2$ are not disjoint, we can use an approach very similar to Section \ref{solution1} in order to derive a bound for Problem \ref{problem2}.

\begin{thm}\label{notdisjoint} Let $F_1$ and $F_2$ be finite subsets of $\mathbb{Z}^2$ such that
\begin{itemize}
\item $F_1$ is uniquely determined by its row and column sums, and
\item $|F_1| = |F_2|$.
\end{itemize}
Let $\alpha$ be defined as in Section \ref{notation}, and let $p = |F_1 \cap F_2|$. Then
\[
|F_1| \leq \sum_{i=1}^{\alpha + p} \left\lfloor \frac{\alpha +p}{i} \right\rfloor.
\]
\end{thm}

\begin{pf}
Assume that the rows and columns are ordered as in Section \ref{notation}. Let $(k,l) \in F_1$. Then all the points in the rectangle $\{ (i,j): 1 \leq i \leq k, 1 \leq j \leq l \}$ are elements of $F_1$. At most $p$ of the points in this rectangle are elements of $F_2$, so at least $kl - p$ points belong to $F_1 \symm F_2$. None of the points in the rectangle is an element of $F_2 \backslash F_1$, so all of the $kl - p$ points of $F_1 \symm F_2$ in the rectangle must belong to different staircases, which implies $\alpha +p \geq kl$. For all $i$ with $1 \leq i \leq a$ we have $(i, r_i^{(1)}) \in F_1$, hence $r_i^{(1)} \leq \frac{\alpha +p}{i}$. Since $r_i^{(1)}$ must be an integer, we have
\[
|F_1| = \sum_{i=1}^a r_i^{(1)} \leq \sum_{i=1}^a \left\lfloor \frac{\alpha +p}{i} \right\rfloor.
\]
Since $(a,1) \in F_1$, we have $a \leq \alpha +p$, so
\[
|F_1| \leq \sum_{i=1}^{\alpha +p} \left\lfloor \frac{\alpha +p}{i} \right\rfloor.
\] \hfill $\square$
\end{pf}

\begin{cor}\label{grenslog}
Let $F_1$, $F_2$, $\alpha$ and $p$ be defined as in Theorem \ref{notdisjoint}. Then
\[
|F_1| \leq (\alpha +p) (1 +\log (\alpha +p)).
\]
\end{cor}

\begin{pf}
Analogous to the proof of Corollary \ref{disjunctegrenslog}. \hfill $\square$
\end{pf}

The following example shows that the upper bound cannot even be improved by a factor $\frac{1}{2 \log 2} \approx 0.72$, provided that $\alpha > \frac{p+1}{2 \log 2 - 1} \log (p+1)$.

\begin{figure}
\begin{center}
\includegraphics{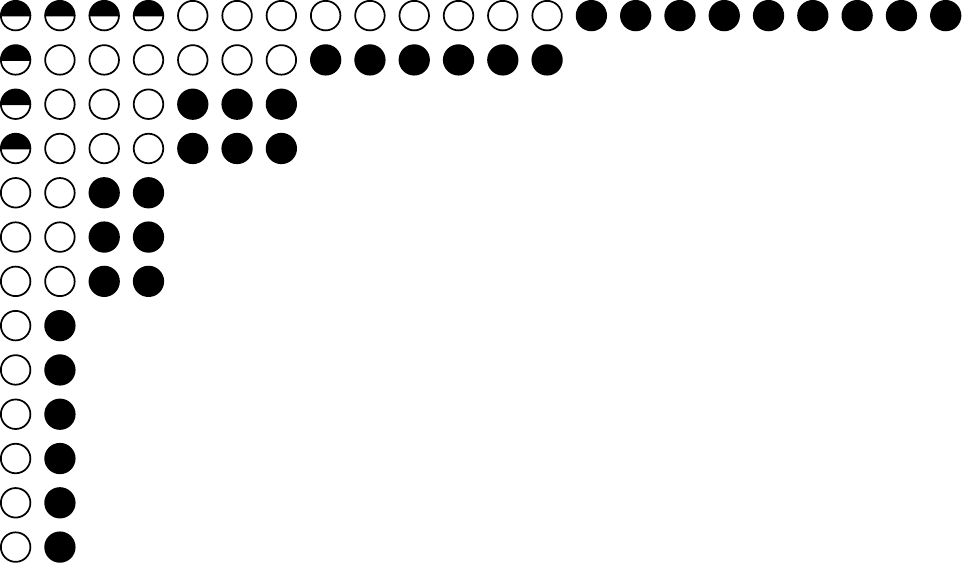}
\end{center}
\caption{The construction from Example \ref{example2} with $k=3$ and $m=4$.}
\label{plaatjevoorbeeld2}
\end{figure}

\begin{figure}
\begin{center}
\includegraphics{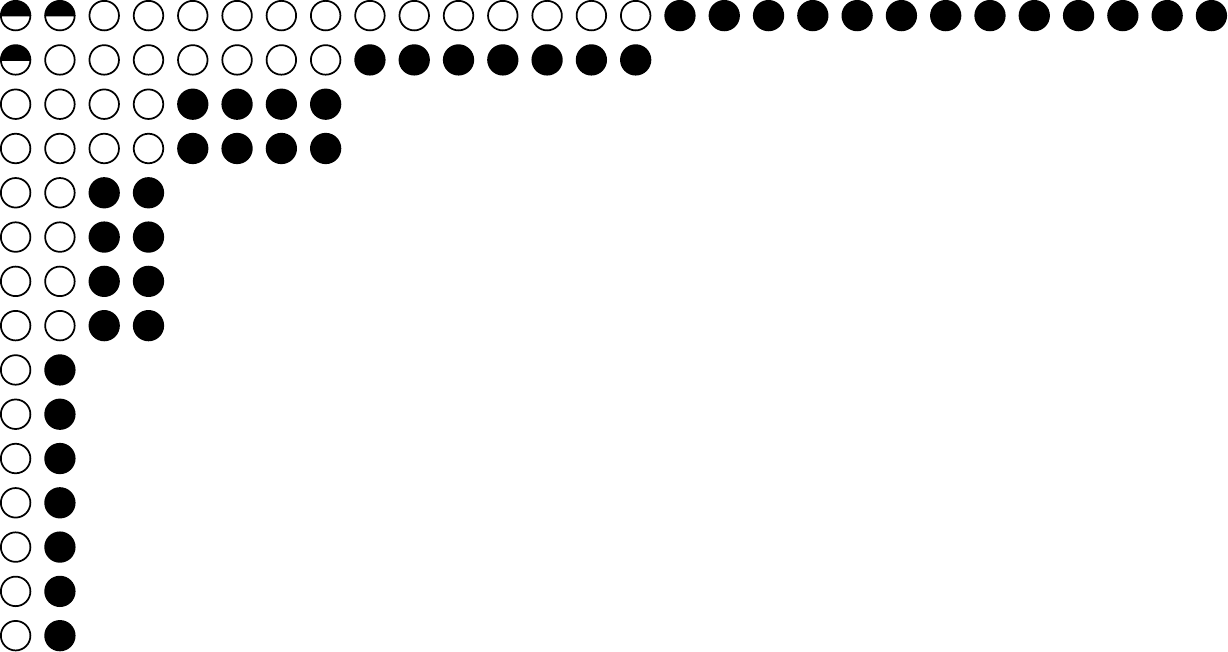}
\end{center}
\caption{The construction from Example \ref{example2} with $k=2$ and $m=4$.}
\label{plaatjevoorbeeld3}
\end{figure}

\begin{exmp}\label{example2}\end{exmp} Let $k$ and $m$ be integers satisfying $k \geq 2$ and $m \geq 2k-2$. We construct sets $F_1$ and $F_2$ as follows (see also Figures \ref{plaatjevoorbeeld2} and \ref{plaatjevoorbeeld3}).
\begin{itemize}
\item Row 1: \begin{itemize} \item $(1,j) \in F_1 \cap F_2$ for $1 \leq j \leq 2^{k-1}$, \item $(1,j) \in F_1$ for $2^{k-1}+1 \leq j \leq 2^m - 2^{k-1} + 1$, \item $(1,j) \in F_2$ for $2^m - 2^{k-1}+2 \leq j \leq 2^{m+1} - 2^k - 2^{k-1} + 2$. \end{itemize}
\item Let $0 \leq l \leq k-2$. Row $i$, where $2^{l}+1 \leq i \leq 2^{l+1}$: \begin{itemize} \item $(i,1) \in F_1 \cap F_2$,  \item $(i,j) \in F_1$ for $2 \leq j \leq 2^{m-l-1} - 2^{k-l-2} +1$, \item $(i,j) \in F_2$ for $2^{m-l-1} - 2^{k-l-2} +2 \leq j \leq 2^{m-l} - 2^{k-l-1} + 1$. \end{itemize}
\item Let $k-1 \leq l \leq m-k$. Row $i$, where $2^l+1 \leq i \leq 2^{l+1}$: \begin{itemize} \item $(i,j) \in F_1$ for  $1 \leq j \leq 2^{m-l-1}$, \item $(i,j) \in F_2$ for $2^{m-l-1}+1 \leq j \leq 2^{m-l}$.  \end{itemize}
\item Let $m-k+1 \leq l \leq m-1$. Row $i$, where $2^l - 2^{l-m+k-1} + 2 \leq i \leq 2^{l+1} - 2^{l-m+k} +1$: \begin{itemize}\item $(i,j) \in F_1$ for $1 \leq j \leq 2^{m-l-1}$, \item $(i,j) \in F_2$ for $2^{m-l-1} +1 \leq j \leq 2^{m-l}$. \end{itemize}
\end{itemize}

The construction is almost symmetrical: if $(i,j) \in F_1$, then $(j,i) \in F_1$; if $(i,j) \in F_1 \cap F_2$, then $(j,i) \in F_1 \cap F_2$; and if $(i,j) \in F_2$ with $i > 1$, then $(j,i) \in F_2$. Since it is clear from the construction that each row contains exactly as many points of $F_1$ as points of $F_2$, we conclude that each column $j$ with $2 \leq j \leq 2^m - 2^{k-1} + 1$ contains exactly as many points of $F_1$ as points of $F_2$ as well. The only difference in the line sums occurs in the first column (which has $2^m - 2^{k-1} + 1$ points of $F_1$ and only  $2^{k-1}$ of $F_2$) and in columns $2^m - 2^{k-1}+2$ up to $2^{m+1} - 2^k - 2^{k-1} + 2$ (each of which contains one point of $F_2$ and none of $F_1$). So we have
\begin{eqnarray*}
\alpha &=& \frac{1}{2} \left( (2^m - 2^{k-1} + 1) - 2^{k-1} + (2^{m+1} - 2^k - 2^{k-1} + 2) - (2^m - 2^{k-1} + 1) \right) \\ &=& 2^m - 2^k + 1.
\end{eqnarray*}
It is easy to see that
\[
p = |F_1 \cap F_2| = 2^k - 1.
\]
Now we count the number of elements of $F_1$.
\begin{itemize}
\item Row 1 contains $2^m - 2^{k-1} + 1$ elements of $F_1$.
\item Let $0 \leq l \leq k-2$. Rows $2^{l}+1$ up to $2^{l+1}$ together contain $2^l (2^{m-l-1} - 2^{k-l-2} +1) = 2^{m-1} - 2^{k-2} + 2^l$ elements of $F_1$.
\item Let $k-1 \leq l \leq m-k$. Rows $2^l+1$ up to $2^{l+1}$ together contain $2^l \cdot 2^{m-l-1} = 2^{m-1}$ elements of $F_1$.
\item  Let $m-k+1 \leq l \leq m-1$. Rows $2^l - 2^{l-m+k-1} + 2$ up to $2^{l+1} - 2^{l-m+k} +1$ together contain $ (2^l - 2^{l-m+k-1})(2^{m-l-1}) = 2^{m-1} - 2^{k-2} $ elements of $F_1$.
\end{itemize}
Hence the number of elements of $F_1$ is
\begin{eqnarray*}
|F_1| &=& 2^m - 2^{k-1} + 1 + (k-1)(2^{m-1} - 2^{k-2}) + \sum_{l=0}^{k-2} 2^l  \\ && \hspace{1cm} + (m-2k+2) 2^{m-1} + (k-1) (2^{m-1} - 2^{k-2}) \\ & = & 2^m + m2^{m-1} + 2^{k-1} - k2^{k-1}.
\end{eqnarray*}
For this family of examples we now have
\[
|F_1| = \alpha + p + \frac{\alpha+p}{2} \log_2(\alpha + p) + \frac{p+1}{2} - \frac{p+1}{2} \log_2(p+1).
\]

We will now prove another bound, which is better if $p = |F_1 \cap F_2|$ is large compared to $\alpha$. Let $u$ be an integer such that $2u = |F_1 \symm F_2|$. We will first derive an upper bound on $u$ in terms of $a$, $b$ and $\alpha$. Then we will derive a lower bound on $|F_1|$ in terms of $a$, $b$ and $\alpha$. By combining these two, we find an upper bound on $u$ in terms of $\alpha$ and $p$.

\begin{lem}\label{vergelijkingm}
Let $F_1$ and $F_2$ be finite subsets of $\mathbb{Z}^2$ such that
\begin{itemize}
\item $F_1$ is uniquely determined by its row and column sums, and
\item $|F_1| = |F_2|$.
\end{itemize}
Let $\alpha$, $a$ and $b$ be defined as in Section \ref{notation}. Define $u$ as $2u = |F_1 \symm F_2|$. Then we have
\[
u^2 \leq \frac{\alpha}{4} (a+b)(a+b+\alpha -1).
\]
\end{lem}

\begin{pf} Decompose $F_1 \symm F_2$ into $\alpha$ staircases as in Lemma \ref{staircases}, and let $\mathcal{T}$ be the set consisting of these staircases. Let $T \in \mathcal{T}$ be a staircase and $i \leq a+1$ a positive integer. Consider the elements of $T \cap F_2$ in rows $i$, $i+1$, \ldots, $a$. If such elements exist, then let $w_i(T)$ be the largest column index that occurs among these elements. If there are no elements of $T \cap F_2$ in those rows, then let $w_i(T)$ be equal to the smallest column index of an element of $T \cap F_1$ (no longer restricted to rows $i$, \ldots, $a$). We have $w_i(T) \geq 1$. Define $W_i = \sum_{T \in \mathcal{T}} w_i(T)$.

Let $d_i$ be the number of elements of $F_1 \backslash F_2$ in row $i$. Let $y_1 < \ldots < y_{d_i}$ be the column indices of the elements of $F_1 \backslash F_2$ in row $i$, and let $y_1' < \ldots < y_{d_i}'$ be the column indices of the elements of $F_2 \backslash F_1$ in row $i$. Let $\mathcal{T}_i \subset \mathcal{T}$ be the set of staircases with elements in row $i$. The elements in $F_2 \backslash F_1$ of these staircases are in columns $y_1'$, $y_2'$, \ldots, $y_{d_i}'$, hence the set $\{w_i(T): T \in \mathcal{T}_i\}$ is equal to the set $\{y_1', y_2', \ldots, y_{d_i}'\}$. The elements in $F_1 \backslash F_2$ are in columns $y_1$, $y_2$, \ldots, $y_d$ and are either the first element of a staircase or correspond to an element of $F_2 \backslash F_1$ in the same column but in a row with index at least $i+1$. In either case, for a staircase $T \in \mathcal{T}_i$ we have $w_{i+1}(T) = y_j$ for some $j$. Hence the set $\{w_{i+1}(T): T \in \mathcal{T}_i\}$ is equal to the set $\{y_1, y_2, \ldots, y_{d_i}\}$. We have
\[
\sum_{T \in \mathcal{T}_i} w_{i+1}(T) = \sum_{j=1}^{d_i} y_j \leq \sum_{j=1}^{d_i} (y_{d_i}-j+1) = d_iy_{d_i} - \frac{1}{2}(d_i-1)d_i,
\]
and
\[
\sum_{T \in \mathcal{T}_i} w_i(T) = \sum_{j=1}^{d_i} y_j' \geq \sum_{j=1}^{d_i} (y_{d_i} + j) = d_iy_{d_i} + \frac{1}{2}(d_i+1)d_i.
\]
Hence
\begin{eqnarray*}
W_i &=& W_{i+1} + \sum_{T \in \mathcal{T}_i} (w_i(T) - w_{i+1}(T)) \\ &\geq & W_{i+1} + \frac{1}{2}(d_i+1)d_i + \frac{1}{2}(d_i-1)d_i \\ &=& W_{i+1} + d_i^2.
\end{eqnarray*}
Since $W_{a+1} \geq \alpha$, we find
\[
W_1 \geq \alpha + d_1^2 + \cdots + d_a^2.
\]

We may assume that if $(x,y)$ is the endpoint of a staircase, then $(x,y')$ is an element of $F_1 \cup F_2$ for $1 \leq y' < y$ (i.e. there are no gaps between the endpoints and other elements of $F_1 \cup F_2$ on the same row). After all, by moving the endpoint of a staircase to another empty position on the same row, the error in the columns can only become smaller (if the new position of the endpoint happens to be in the same column as the first point of another staircase, in which case the two staircases fuse together to one) but not larger, and $u$, $a$ and $b$ do not change.

So on the other hand, as $W_1$ is the sum of the column indices of the endpoints of the staircases, we have
\[
W_1 \leq (b+1) + (b+2) + \cdots + (b+\alpha) = \alpha b + \frac{1}{2}\alpha(\alpha+1).
\]
We conclude
\[
\alpha + \sum_{i=1}^a d_i^2 \leq \alpha b + \frac{1}{2}\alpha(\alpha+1).
\]
Note that $\sum_{i=1}^a d_i = u$. By the Cauchy-Schwarz inequality, we have
\[
\left( \sum_{i=1}^a d_i^2 \right) \left( \sum_{i=1}^a 1 \right) \geq \left( \sum_{i=1}^a d_i \right)^2 = u^2,
\]
so
\[
\sum_{i=1}^a d_i^2 \geq \frac{u^2}{a}.
\]
From this it follows that
\[
\alpha b + \frac{1}{2}\alpha(\alpha+1) \geq \alpha + \frac{u^2}{a},
\]
or, equivalently,
\[
u^2 \leq \alpha ab + \frac{1}{2}\alpha(\alpha-1)a.
\]
By symmetry we also have
\[
u^2 \leq \alpha ab + \frac{1}{2}\alpha(\alpha-1)b.
\]
Hence
\[
u^2 \leq \alpha ab + \frac{1}{4}\alpha(\alpha-1)(a+b).
\]
Using that $\sqrt{ab} \leq \frac{a+b}{2}$, we find
\[
u^2 \leq \alpha \left( \frac{(a+b)^2}{4} + \frac{(\alpha-1)(a+b)}{4}  \right) = \frac{\alpha}{4}(a+b)(a+b+\alpha-1).
\] \hfill $\square$
\end{pf}

\begin{lem}\label{vergelijkingF}
Let $F_1$ and $F_2$ be finite subsets of $\mathbb{Z}^2$ such that
\begin{itemize}
\item $F_1$ is uniquely determined by its row and column sums, and
\item $|F_1| = |F_2|$.
\end{itemize}
Let $\alpha$, $a$ and $b$ be defined as in Section \ref{notation}. Then
we have
\[
|F_1| \geq \frac{(a+b)^2}{4(\alpha+1)}.
\]
\end{lem}

\begin{pf}
Without loss of generality, we may assume that all rows and columns that contain elements of $F_1$ also contain at least one point $F_1 \symm F_2$: if a row or column does not contain any points of $F_1 \symm F_2$, we may delete it. By doing so, $F_1 \symm F_2$ does not change, while $|F_1|$ becomes smaller, so the situation becomes better.

First consider the case $r_{i+1}^{(1)} < r_i^{(1)} - \alpha$ for some $i$. We will show that this is impossible. If a column does not contain an element of $F_2 \backslash F_1$, then by the assumption above it contains an element of $F_1 \backslash F_2$, which must then be the first point of a staircase. Consider all points of $F_2 \backslash F_1$ and all first points of staircases in columns $r_{i+1}+1$, $r_{i+1}+2$, \ldots, $r_i$. Since these are more than $\alpha$ columns, at least two of those points must belong to the same staircase. On the other hand, if $(x,y) \in F_1 \backslash F_2$ is the first point of a staircase with $r_{i+1} < y \leq r_i$, then we have $x \leq i$, so the second point $(x', y')$ in the staircase, which is in $F_2 \backslash F_1$, must satisfy $x' \leq i$ and therefore $y' > r_i$. So the second point cannot also be in one of the columns $r_{i+1}+1$, $r_{i+1}+2$, \ldots, $r_i$. If two points of $F_2 \backslash F_1$ in columns $r_{i+1}+1$, $r_{i+1}+2$, \ldots, $r_i$ belong to the same staircase, then they must be connected by a point of $F_1 \backslash F_2$ in the same columns. However, by a similar argument this forces the next point to be outside the mentioned columns, while we assumed that it was in those columns. We conclude that it is impossible for row sums of two consecutive rows to differ by more than $\alpha$.

By the same argument, column sums of two consecutive columns cannot differ by more than $\alpha$. Hence we have $r_{i+1}^{(1)} \geq r_i^{(1)} - \alpha$ for all $i$, and $c_{j+1}^{(1)} \geq c_j^{(1)} - \alpha$ for all $j$.

We now have $r_2^{(1)} \geq b - \alpha$, $r_3^{(1)} \geq b - 2 \alpha$, and so on. Also, $c_2^{(1)} \geq a - \alpha$, $c_3^{(1)} \geq a - 2\alpha$, and so on. Using this, we can derive a lower bound on $|F_1|$ for fixed $a$ and $b$. Consider Figure \ref{plaatjeoppervlakte}. The points of $F_1$ are indicated by black dots. The number of points is equal to the grey area in the picture, which consists of all $1 \times 1$-squares with a point of $F_1$ in the upper left corner. We can estimate this area from below by drawing a line with slope $\alpha$ through the point $(a+1, 1)$ and a line with slope $\frac{1}{\alpha}$ through the point $(b+1,1)$; the area closed in by these two lines and the two axes is less than or equal to the number of points of $F_1$.

\begin{figure}
\begin{center}
\includegraphics{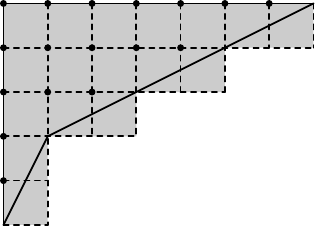}
\end{center}
\caption{The number of points of $F_1$ (indicated by small black dots) is equal to the grey area.}
\label{plaatjeoppervlakte}
\end{figure}

For $\alpha =1$ those lines do not have a point of intersection. Under the assumption we made at the beginning of this proof, we must in this case have $a=b$ and the number of points of $F_1$ is equal to
\[
\frac{a(a+1)}{2} \geq \frac{a^2}{\alpha+1} = \frac{(a+b)^2}{4(\alpha+1)},
\]
so in this case we are done.

In order to compute the area for $\alpha \geq 2$ we switch to the usual coordinates in $\mathbb{R}^2$, see Figure \ref{plaatjeoppervlakte2}. The equation of the first line is $y= \alpha x - a$, and the equation of the second line is $y= \frac{1}{\alpha} x - \frac{1}{\alpha} b$. We find that the point of intersection is given by
\[
(x,y) = \left( \frac{a\alpha - b}{\alpha^2 -1} , \frac{-b\alpha +a }{\alpha^2 -1} \right).
\]
The area of the grey part of Figure \ref{plaatjeoppervlakte2} is equal to
\[
\frac{1}{2} a \cdot \frac{a \alpha - b}{\alpha^2-1}\ +\ \frac{1}{2} b \cdot \frac{b\alpha - a}{\alpha^2 -1} \ =\ \frac{ a^2\alpha + b^2 \alpha - 2ab }{2(\alpha^2-1)}.
\]
We now have
\[
|F_1| \geq \frac{\alpha( a^2 + b^2) - 2ab }{2(\alpha^2-1)} \geq \frac{\alpha \frac{(a+b)^2}{2} - \frac{(a+b)^2}{2}}{2(\alpha^2-1)} = \frac{(a+b)^2}{4(\alpha+1)}.
\]
\hfill $\square$

\begin{figure}
\begin{center}
\includegraphics{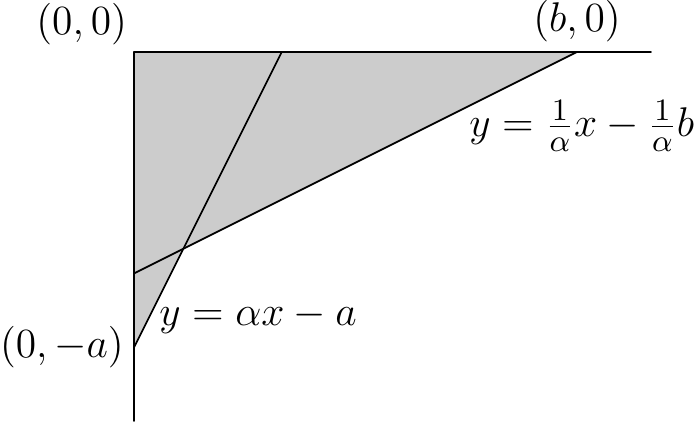}
\end{center}
\caption{Computing the area bounded by the two lines and the two axes.}
\label{plaatjeoppervlakte2}
\end{figure}

\end{pf}

\begin{thm}\label{grenswortel}
Let $F_1$ and $F_2$ be finite subsets of $\mathbb{Z}^2$ such that
\begin{itemize}
\item $F_1$ is uniquely determined by its row and column sums, and
\item $|F_1| = |F_2|$.
\end{itemize}
Let $\alpha$ be defined as in Section \ref{notation}, and let $p = |F_1 \cap F_2|$. Write $\beta = \sqrt{\alpha}(\alpha+1)$. Then
\[
|F_1| \leq p + \sqrt{\frac{\alpha}{4} \left(
 \beta + \sqrt{\beta(\alpha-1) + 4(\alpha+1)p + \beta^2  } + \frac{\alpha-1}{2} \right)^2 - \frac{(\alpha-1)^2 \alpha}{16} }.
\]
\end{thm}

\begin{pf}
Write $s = a+b$ for convenience of notation. From Lemma \ref{vergelijkingm} we derive
\[
u \leq \frac{\sqrt{\alpha}}{2} \left( s + \frac{\alpha-1}{2} \right).
\]
We substitute $|F_1| = u+p$ in  Lemma \ref{vergelijkingF} and use the above bound for $u$:
\[
\frac{\sqrt{\alpha}}{2} \left( s + \frac{\alpha-1}{2} \right) + p \geq |F_1| \geq  \frac{s^2}{4(\alpha+1)}.
\]
Solving for $s$, we find
\begin{eqnarray*}
s & \leq & \sqrt{\alpha}(\alpha+1) + \sqrt{\sqrt{\alpha}(\alpha^2-1) + 4(\alpha+1)p + \alpha(\alpha+1)^2  } \\
& = & \beta + \sqrt{\beta(\alpha+1) + 4(\alpha+1)p + \beta^2  }
\end{eqnarray*}
Finally we substitute this in Lemma \ref{vergelijkingm}:
\[
u \leq \sqrt{\frac{\alpha}{4} \left(
\beta + \sqrt{\beta(\alpha-1) + 4(\alpha+1)p + \beta^2  } + \frac{\alpha-1}{2} \right)^2 - \frac{(\alpha-1)^2 \alpha}{16} }.
\]
This, together with $|F_1| = u+p$, yields the claimed result.
\hfill $\square$
\end{pf}

\begin{rem}
By a straightforward generalisation of \cite[Proposition 13 and Lemma 16]{alpersartikel}, we find a bound very similar to the one in Theorem \ref{grenswortel}:
\[
|F_1| \leq p+ (\alpha + 1)(\alpha - \frac{1}{2}) + (\alpha + 1) \sqrt{2p + \frac{(2\alpha-1)^2}{4}}.
\]
\end{rem}

Theorem \ref{grenswortel} says that $|F_1|$ is asymptotically bounded by $p + \alpha \sqrt{p} + \alpha^2$. The next example shows that $|F_1|$ can be asymptotically as large as $p + 2\sqrt{\alpha p} + \alpha$.

\begin{figure}
\begin{center}
\includegraphics{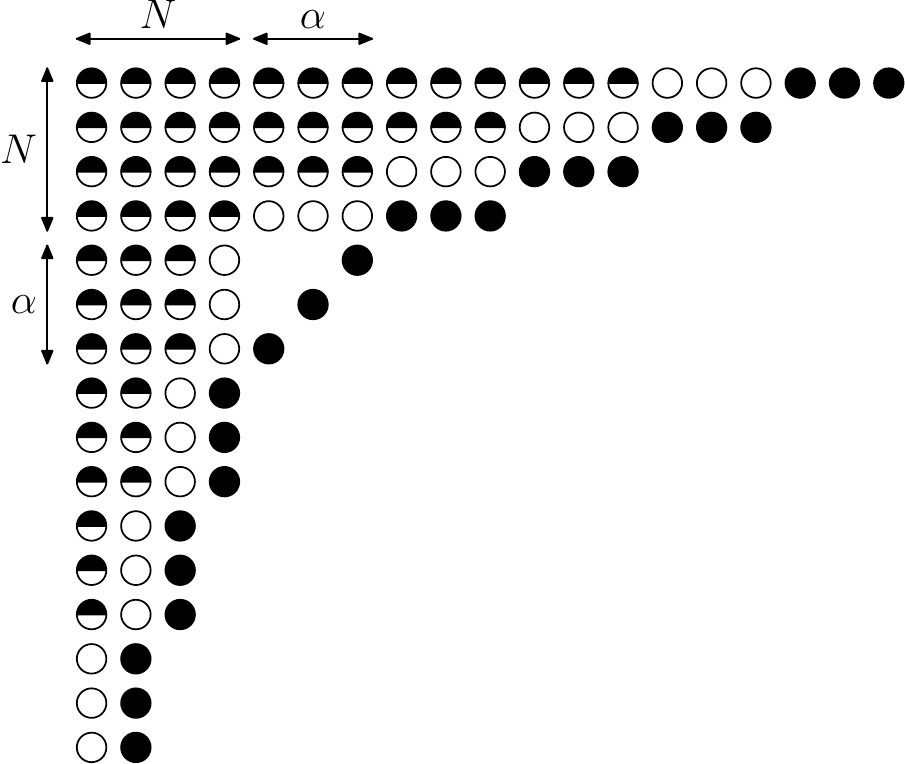}
\end{center}
\caption{The construction from Example \ref{example3} with $N=4$ and $\alpha =3$.}
\label{plaatjevoorbeeld4}
\end{figure}

\begin{exmp}\label{example3} \end{exmp} Let $N$ be a positive integer. We construct $F_1$ and $F_2$ with total difference in the line sums equal to $2 \alpha$ as follows (see also Figure \ref{plaatjevoorbeeld4}).
Let $(i,j) \in F_1 \cap F_2$ for $1 \leq i \leq N$, $1 \leq j \leq N$. Furthermore, for $1 \leq i \leq N$:
\begin{itemize}
\item Let $(i,j), (j,i) \in F_1 \cap F_2$ for $N+1 \leq j \leq N+(N-i)\alpha$.
\item Let $(i,j), (j,i) \in F_1$ for $N+(N-i)\alpha + 1 \leq j \leq N+(N-i+1)\alpha$.
\item Let $(i,j), (j,i) \in F_2$ for $N+(N-i+1)\alpha +1 \leq j \leq N+(N-i+2)\alpha$.
\end{itemize}
Finally, for $1 \leq t \leq \alpha$, let $(i,j) \in F_2$ with $i=N+t$ and $j = N+\alpha + 1 -t$.
\\

The only differences in the line sums occur in the first column (a difference of $\alpha$) and in columns $N+N\alpha + 1$ up to $N +N\alpha + \alpha$ (a difference of 1 in each column). We have
\[
p = N^2 + 2 \cdot \frac{1}{2}N(N-1)\alpha = N^2 + N^2 \alpha - N \alpha,
\]
and
\[
|F_1| = N^2 + 2 \cdot \frac{1}{2}N(N+1)\alpha = N^2 + N^2 \alpha + N\alpha.
\]
From the first equality we derive
\[
N = \frac{\alpha}{2(\alpha+1)} + \sqrt{\frac{p}{\alpha+1} + \frac{\alpha^2}{4(\alpha+1)^2}}.
\]
Hence
\[
|F_1| = p + 2N\alpha = p + \frac{\alpha^2}{\alpha+1} + \sqrt{\frac{4\alpha^2 p}{\alpha+1} + \frac{\alpha^4}{(\alpha+1)^2}}.
\]

\section{Generalisation to unequal sizes}\label{neq}

Until now, we have assumed that $|F_1| = |F_2|$. However, we can easily generalise all the results to the case $|F_1| \neq |F_2|$.

Suppose $|F_1| > |F_2|$. Then there must be a row $i$ with $r_i^{(1)} > r_i^{(2)}$. Let $j > b$ be such that $(i,j) \not\in F_2$ and define $F_3 = F_2 \cup \{(i,j)\}$. We have $r_i^{(3)} = r_i^{(2)} + 1$, so the error in row $i$ has decreased by one, while the error in column $j$ has increased by one. In this way, we can keep adding points until $F_2$ together with the extra points is just as large as $F_1$, while the total difference in the line sums is still $2 \alpha$. Note that $p = |F_1 \cap F_2|$ and $|F_1|$ have not changed during this process, so the results from Theorem \ref{grenswortel} and Corollary \ref{grenslog} are still valid in exactly the same form.

Suppose on the other hand that $|F_1| < |F_2|$. Then there must be a row with $r_i^{(1)} < r_i^{(2)}$. Let $j$ be such that $(i,j) \in F_2 \backslash F_1$ and define $F_3 = F_2 \backslash \{(i,j)\}$. The error in row $i$ has now decreased by one, while the error in column $j$ has at most increased by one, so the total error in the line sums has not increased. We can keep deleting points of $F_2$ until there are exactly $|F_1|$ points left, while the total difference in the line sums is at most $2\alpha$.

By using $|F_1 \symm F_2| = 2(|F_1| - p)$, we can state the results from Theorem \ref{grenswortel} and Corollary \ref{grenslog} in a more symmetric way, not depending on the size of $F_1$.

\begin{thm}
Let $F_1$ and $F_2$ be finite subsets of $\mathbb{Z}^2$ such that
$F_1$ is uniquely determined by its row and column sums.
Let $\alpha$ be defined as in Section \ref{notation}, and let $p = |F_1 \cap F_2|$. Write $\beta = \sqrt{\alpha}(\alpha+1)$. Then
\begin{enumerate}
\item $|F_1 \symm F_2| \leq 2\alpha + 2(\alpha+p) \log(\alpha+p).$
\item $|F_1 \symm F_2| \leq \sqrt{\alpha \left(
\beta + \sqrt{\beta(\alpha-1) + 4(\alpha+1)p + \beta^2  } + \frac{\alpha-1}{2} \right)^2 - \frac{(\alpha-1)^2 \alpha}{4} }.$
\end{enumerate}
\end{thm}

\end{document}